\documentclass[12pt]{amsart}
\usepackage{amsmath,amscd,amsthm,amsfonts, amssymb,amsxtra}
\usepackage{epsf}
\usepackage[all]{xy}
\newtheorem{thm}{Theorem}[section]  
\newtheorem*{un-no-thm}{Theorem}
\newtheorem{cor}[thm]{Corollary}     % Numbered along with thm
\newtheorem{lem}[thm]{Lemma}         % Numbered along with thm
  
\newtheorem{bigthm}{Theorem}
\newtheorem{bigcor}[bigthm]{Corollary}

\theoremstyle{definition}
\theoremstyle{definition}
\theoremstyle{definition}

\theoremstyle{remark}

\newtheorem*{acks}{Acknowledgements}
 
\newtheorem*{out}{Outline}

\newtheorem*{rems_no}{Remarks}
\newtheorem*{rem_no}{Remark}

\newtheorem*{convent}{Conventions}

\begin{document}
\title{On embeddings in the sphere}
\date{\today}
\author{John R. Klein}
\address{Wayne State University, Detroit, MI 48202}
\email{klein@math.wayne.edu}
\begin{abstract} We consider 
embeddings of a finite complex 
in a sphere. 
We give a homotopy theoretic classification such embeddings in a wide range.
\end{abstract}
\thanks{The author is partially supported by  NSF
Grant DMS-0201695.}
\maketitle%\pagestyle{fancy}
\setlength{\parindent}{15pt}
\setlength{\parskip}{1pt plus 0pt minus 1pt}
\def\Top{\bold T\bold o \bold p}
\def\Sp{\bold S\bold p}
\def\vo{\varOmega}
\def\vs{\varSigma}
\def\smsh{\wedge}
\def\flush{\flushpar}
\def\id{\text{id}}
\def\dbslash{/\!\! /}
\def\codim{\text{\rm codim\,}}
\def\:{\colon}
\def\holim{\text{\rm holim\,}}
\def\hocolim{\text{\rm hocolim\,}}
\def\hodim{\text{\rm hodim\,}}
\def\hocodim{\text{hocodim\,}}
\def\Bbb{\mathbb}
\def\bold{\mathbf}
\def\Aut{\text{\rm Aut}}
\def\cal{\mathcal}
\def\frak{\mathfrak}

\section{Introduction}
Let $K$ be a finite complex.  An
{\it  embedding up to homotopy} of $K$ in $S^n$
consists of a a pair
$$
(M,h)
$$
where $M^n$ is compact 
codimension zero PL submanifold of $S^n$ and 
$h\:K \to M$ is a homotopy equivalence. Two 
such pairs $(M_0,h_0)$ and $(M_1,h_1)$
are said to be {\it concordant} if there is an embedded $h$-cobordism $W$ in
$S^n{\times}[0,1]$ from $M_0$ to $M_1$ together with a homotopy equivalence
$H\:K{\times}[0,1]\to W$ extending both $h_0$ and $h_1$.
Let 
$$
E(K,S^n)
$$ 
denote the set of concordance classes of embeddings up to homotopy
of $K$ in $S^n$. (Note: if $K$ is $1$-connected and $\dim K \le n{-}3$,
the existence of a concordance implies the existence of an ambient isotopy.)
Unless confusion arises, we refer to embeddings up to homotopy as {\it embeddings.}

\subsection*{Constraints} We fix throughout integers 
$$
k, n, r
$$ 
satisfying
$$
0 \le k \le n{-}3, \quad r \ge 1, \quad \text{ and } \quad n\ge 6\, .
$$
If
$n \le 7$, we also assume $k - r \ge 2$.
 
In addition to these constraints, we 
consider the inequalities
\begin{equation}\label{number1}
r \ge \max (\frac{1}{2}(2k{-}n),3k {-}2n{+}2) 
\end{equation}

\begin{equation} \label{number2}
r \ge \max (\frac{1}{2}(2k{-}n{+}1),3k {-}2n{+}3) 
\end{equation}

The inequalities can be interpreted as 
follows: the integer $r$ will be the connectivity 
of the space to be embedded. 
Consider maps from manifolds of dimension $k$ to $S^n$.
Then roughly, the inequalities represent 
the demand that the connectivity $r$ exceeds
both the generic dimension of the triple point set and also,
one half the generic dimension of the double point set.

\subsection*{Main results} To formulate our main results
requires some preparation.
Let $$\alpha\:\Bbb Z_2 \to {\rm GL_1(\Bbb R)}$$ denote the 
sign representation. If $s,t \ge 0$ denote integers, let $S^{t\alpha + s}$
denote the one point compactification of the direct sum of $t$ copies of
$\alpha$ with $s$-copies of the
trivial representation. 
This is a sphere of dimension $t+s$ having a based action of 
$\Bbb Z_2$. 

If $X$ and $Y$ are based spaces, we let 
$F^{\rm st}(X,Y)$ denote the spectrum
of stable maps from $X$ to $Y$ (the $j$-th space
of this spectrum is the function space of maps
from $X$ to $Q(\Sigma^j Y)$, where $Q$ denotes the
stable homotopy functor).

If $X$ and $Y$ are based  $\Bbb Z_2$-spaces,
then $F^{\rm st}(X,Y)$ comes equipped the structure of a naive 
$\Bbb Z_2$-spectrum by conjugating functions. 
Let $F^{\rm st}(X,Y)_{h\Bbb Z_2}$ denote the
associated homotopy orbit spectrum. 

Choose a basepoint for $K$. We consider the case
when $X = K\smsh K$ with permutation action and 
$Y = S^{(n{-}1)\alpha + 1}\smsh K$
with the diagonal action (where $\Bbb Z_2$ acts trivially on
$K$). 

We are now in a position to state our main results.

\begin{bigthm}[Existence]\label{exist} 
Assume $K$ is $r$-connected and $\dim K \le k$.
There is an obstruction 
$$
\theta_K \in \pi_0(F^{\rm st}(K\smsh K,S^{(n{-}1)\alpha + 1}\smsh K)_{h\Bbb Z_2})\, ,
$$
depending only on the homotopy type of $K$,
whose vanishing is a necessary condition for $E(K,S^n)$ 
to be non-empty.

If the inequality \eqref{number1} holds, then the 
vanishing of $\theta_K$ implies $E(K,S^n)$ is non-empty.
\end{bigthm}
\medskip

\begin{rems_no} When $K$ is $(2k{-}n)$-connected, the obstruction 
group is trivial, so there is an embedding of $K$ in $S^n$.
Thus we recover the Stallings-Wall embedding theorem (see
\cite{Wall4}, \cite{Stallings}).
\medskip

When $K$ is $(2k{-}n{-}1)$-connected, the
obstruction group is isomorphic to
$$
H^{2k}(K {\times} K;\pi_{2k-n}(K))/(1-T)\, ,
$$
where $T$ is the involution on $H^{2k}(K {\times} K;\pi_{2k-n}(K))$
given by $t\circ E$, where $E$  switches the factors of $K\times K$, and
$t$ is the involution of $\pi_{2k-n}(K)$ given by multiplication by
$(-1)^{n{-}1}$.

This abelian group appears in the work of Habegger \cite{Habegger_embed}, who
gave necessary and sufficient conditions for finding 
embeddings in the fringe dimension beyond the Stallings-Wall
range. Habegger 
defined his obstruction using PL intersection theory.
By contrast, our result will be derived homotopy theoretically by
from an embedding theorem of Connolly and Williams \cite{Connolly-Williams}.
\end{rems_no} 
\bigskip

\begin{bigthm}[Enumeration]\label{enumerate} 
Let $K$ be as above. Fix a basepoint in $E(K,S^n)$.
Then there is a function
$$
\phi_K\: E(K,S^n) \to \pi_0(F^{\rm st}(K\smsh K,S^{(n{-}1)\alpha}\smsh K)_{h\Bbb Z_2})
$$
which is onto if  inequality \eqref{number1} holds.
If inequality \eqref{number2} holds,
then $\phi_K$ is also one to one.
\end{bigthm}
\medskip

\begin{bigcor}[Group Structure] Assume $E(K,S^n)$ is 
 non-empty. If inequality \eqref{number2} holds, then  
$E(K,S^n)$ has the structure of an abelian group.
\end{bigcor}
\bigskip

The above results have corollaries which are too
numerous to list in this introduction
(see \S5-7). For example, here is
a consequence of Theorem \ref{enumerate} which appears
to be new (cf.\ \ref{pd_embed_finiteness} below):

\begin{bigcor}[Isotopy Finiteness] \label{isotopy}
In the range of inequality \eqref{number2},
an $r$-connected closed PL manifold $M^k$ with
trivial betti number $b_{2k{-}n{+}1}(M)$ admits
only finitely many locally flat embeddings in $S^n$ up to isotopy.
\end{bigcor}
\medskip

\begin{out} In \S2 we recall the statement of the Connolly-Williams
Classification Theorem. In \S3 we prove Theorem \ref{enumerate}.
In \S4 we prove Theorem \ref{exist} 
by modifying the proof of Theorem \ref{enumerate}. \S5 contains
applications to embeddings of complexes with 2-cells (these
applications are already in the literature in some form). In \S6
we give applications to embeddings of Poincar\'e spaces and manifolds
(many of the results in this section are new to the literature). 
In \S7 we show that
the obstructions to embedding in the range of inequality
\eqref{number1} are $2$-local.
\end{out}
\bigskip

\begin{convent} We work within the category of compactly
generated (based) spaces. Products are to be re-topologized using
the compactly generated topology.
A space is {\it homotopy finite} if it
is the retract of a finite cell complex.

A non-empty space $X$ is {\it $r$-connected} if 
its homotopy vanishes in degrees $\le r$ for every
choice of basepoint. Note that every non-empty space is 
$({-}1)$-connected. By convention, the empty space is $({-}2)$-connected. 

A map $X \to Y$ (with $Y$ non-empty)
is $r$-connected if it's homotopy fiber at every 
choice of basepoint is $(r{-}1)$-connected.
A {\it weak equivalence} is a map which is $r$-connected
for every $r$.

We write $\dim X \le n$ if $X$ is weak equivalent
to a cell complex having cells of dimension at most $n$.
\end{convent}
\bigskip

\begin{acks} I wish to thank Bill Richter for
introducing me to the notion of Poincar\'e embedding.
Bill also gave me a copy of Habegger's thesis
to read when I was an undergraduate in the early 1980s.
I am very much indebted to Bruce Williams for introducing
me to his papers on embeddings. I am also 
grateful to Bruce for
the wealth of mathematical discussion that I've shared with him.
\end{acks}

\section{The Connolly-Williams classification theorem}

We unearth a
result of Connolly and Williams
which relates $E(K,S^n)$ to a desuspension question.

For a $1$-connected homotopy finite space $K$, consider the set of pairs
$(C,\alpha)$ where $C$ is a $1$-connected homotopy finite space and
$$
\alpha\: S^n \to K\ast C
$$ (the join) induces, via the slant
product, an isomorphism in reduced singular (co-)homology
$\tilde H^*(K) \cong \tilde H_{n{-}*{-}1}(C)$. Introduce
an equivalence relation on such pairs by declaring that
$(C,\alpha) \sim (C',\alpha')$ if and only if there is 
a homotopy equivalence of spaces $g\:C \to C'$
satisfying $(\id_K \ast g)\circ \alpha \simeq \alpha'$. 
Call the resulting set of equivalence classes 
$SW_n(K)$.

There is an evident map of sets
$$
E(K,S^n) \to SW_n(K)
$$
which assigns to an embedding $(M,h)$ of $K$ 
the complement of a choice regular neighborhood of $M$
together its Spanier-Whitehead duality pairing.

\begin{thm}[Connolly-Williams \cite{Connolly-Williams}] \label{frank-bruce} Assume that
$K$ is $r$-connected ($r \ge 1$) and  $\dim K \le k$.
Furthermore, assume $k \le n{-}3$, $n \ge 6$
and $2(k-r) \le n$; if $n \le 7$ assume $k - r \le 2$.
Then 
$$
E(K,S^n) \to SW_n(K)
$$
is onto. If in addition, $2(k-r) \le n-1$. The map is one to one.
\end{thm}
\medskip

\begin{rems_no} On the face of it, this result doesn't provide a
``classification'' of embeddings. Indeed, it isn't clear whether
$SW_n(K)$ is non-empty. The remainder of this paper
 will be concerned with the problem of determining $SW_n(K)$ 
when additional constraints are present.
\medskip

The Connolly-Williams result requires $n \ge 6$ 
because surgery theory is used in the proof. 
A Poincar\'e embedding
version of this result also holds without the requirement $\ge 6$
or additional conditions in dimensions $\le 7$.
The Poincar\'e version can be proved with the
fiberwise homotopy theoretic techniques 
appearing in \cite{Klein_haef}. 
I intend to give a proof of the Poincar\'e
version in a future paper.
\end{rems_no}
\medskip

\subsection*{A variant} We next describe a variant of
$SW_n(K)$ which is more convenient to work with.
Assume that $K$ is equipped with a basepoint.

Let $\cal D_{n{-}1}(K)$ be defined as follows: consider the
set of pairs $(W,\alpha)$ such that $W$ is a based space
and $\alpha\:S^{n-1} \to K \smsh W$ is a stable $S$-duality map.
Define equivalence relation by $(W,\alpha) \sim (W',\alpha')$
if and only if there is an (unstable, based)  map
$g\:W \to W'$ such that $(\id_K \smsh g) \circ \alpha \simeq \alpha'$.

\begin{lem} \label{sw=d} Assume that $K$ is $r$-connected
($r \ge 1$), $\dim K \le k$ and $k \le n{-}3$. Then there is a function
$$
\phi \:SW_n(K) \to \cal D_{n{-}1}(K)
$$ 
which is onto if $2(k{-}r)\le n + 1$. If $2(k{-}r)\le n$, then
$\phi$ is also one to one.
\end{lem}

\begin{proof} 
Let $(C,\alpha)$ be a representative of $SW_n(K)$.
Choose a basepoint for $C$. There is a well-known natural
weak equivalence 
$$
K \ast C \overset\sim\to \Sigma K \smsh C \, . 
$$
Precomposing this weak equivalence with the
map $\alpha$, we obtain a map $S^n \to \Sigma K \smsh C$
which we can arrange to be a based map by precomposing
with a suitable rotation. The associated stable
map $S^{n{-}1} \to K \smsh C$ is an $S$-duality.
We leave it to the reader to check that $\phi$ is well-defined.

We now check that $\phi$ is onto. Let $(W,\alpha)$ respresent
an element of $\cal D_{n-1}(K)$. Then $\alpha\: S^{n{-}1} \to
K \smsh W$ is a stable $S$-duality map. It follows that
$\tilde H_*(W) \cong \tilde H^{n{-}*{-}1}(K) = 0$ if $n{-}*{-}1 > k$.
Thus $W$ has vanishing homology when $* \le  n{-}k{-}2$.
 In particular, as $k \le n{-}3$, it follows that $H_1(W) = 0$.

Let $i\:W \to W^+$ natural map to the (Quillen) plus construction. 
Then $W^+$ is $1$-connected and we have
$$
(W,\alpha) \sim (W^+,\id_K \smsh i) \, .
$$
Using $S$-duality, it is also
straightforward to check that $W^+$ is homotopy finite.
Consequently, we are entitled to assume without loss in generality that
$W$ is $1$-connected and homotopy finite. 

In fact, the above argument shows that
$W$ is $(n{-}k{-}2)$-connected. We infer that the smash
product $\Sigma K \smsh W$ is $(n{-}k{+}r)$-connected.
By the Freudenthal suspension theorem, the stable
map $S^{n{-}1} \to K \smsh W$ is represented 
by an unstable map $\beta\: S^n \to \Sigma K \smsh W$ when $2(k{-}r)\le n + 1$
(unique up to homotopy if $2(k{-}r)\le n$).
This shows that the function $\phi$ is onto if $2(k{-}r)\le n + 1$.
This argument also shows that $\phi$ is one to one if  
$2(k{-}r)\le n$.
\end{proof}

\begin{cor} \label{reformulation}
The statement of Theorem \ref{frank-bruce}
holds when $SW_n(K)$ is replaced by $\cal D_{n{-}1}(K)$.
\end{cor}

\section{Proof of Theorem \ref{enumerate}}

Theorem \ref{enumerate}
will follow from an enumeration result for suspension spectra
appearing in \cite{Klein_susp}. We first review
the statement of this result.

Fix a $1$-connected spectrum $E$. For technical
reasons, we shall assume that $E$ is an $\Omega$-spectrum,
and that spaces of the spectrum $E_j$ are cofibrant
(i.e., retracts of cell complexes).
Consider the set of pairs $$(X,h)$$ such that $X$ is 
a based space
and $h \: \Sigma^{\infty} X \to E$ is a weak 
(homotopy) equivalence.
Define
$$
(X,h)  \sim (Y,g)
$$ if there is a map
of spaces $f\:X \to Y$ such that $g\circ \Sigma^{\infty}f$ is homotopic
to $h$ (in particular, $f$ is
a homology isomorphism). This generates an equivalence relation.
Let $\Theta_E$ denote the associated set of equivalence classes.

We write $\dim E \le k$ if $E$ can be obtained from the
trivial spectrum by attaching cells of dimension $\le k$.
Recall that the {\it second extended power}
$D_2(E)$ is the homotopy orbit spectrum of $\Bbb Z_2$ acting
on $E^{\smsh 2}$.

\begin{thm}[Klein \cite{Klein_susp}]\label{desuspend}
Assume $\Theta_E$ is nonempty
and is equipped with a choice of basepoint. 
Then there is a basepoint preserving function
$$
\phi\:\Theta_E \to [E,D_2(E)]\, .
$$
If $E$ is $r$-connected, $r\ge 1$ and  $\dim E \le 3r+2$,
Then $\phi$ is a surjection. If in addition $\dim E \le 3r+1$,
$\phi$ is a bijection.
\end{thm}

\subsection{} 
Recall that 
$$
F^{\rm st}(K, S^{n{-}1})
$$
is spectrum of stable maps from $K$ to $S^{n{-}1}$. 

\begin{lem} \label{first} There is a bijection
$$
\Theta_{F^{\rm st}(K,S^{n{-}1})} \cong \cal D_{n-1}(K)\, .
$$
\end{lem} 

\begin{proof} An element of $\Theta_{F^{\rm st}(K,S^{n{-}1})}$ 
is represented by a pair $(C,\alpha)$, where
$C$ is a based space and  a weak equivalence $\alpha \: \Sigma^{\infty} C 
\to F^{\rm st}(K,S^{n{-}1})$. Taking the adjunction, this is
the same as specifying a (stable) $S$-duality map $\alpha \:
C \smsh K \to S^{n-1}$.
A standard application of $S$-duality (the ``umkehr'' or transpose
map) then allows us
to associate to $\alpha$ an $S$-duality map $\alpha^*\: S^{n-1}\to K \smsh C$. 
The pair $(C,\alpha^*)$ then represents an element of $\cal D_{n{-}1}(K)$.
It is straightforward to check that this procedure defines a bijection.
\end{proof}

\begin{lem} \label{second} Let $E = F^{\rm st}(K,S^{n-1})$. Then
there is an isomorphism of abelian groups
$$
[E,D_2(E)] \cong 
\pi_0(F^{\rm st}(K\smsh K,S^{(n{-}1)\alpha}\smsh K)_{h\Bbb Z_2})
$$
\end{lem} 

\begin{proof} It will be convenient for
 us to rewrite $E \simeq  K^*\smsh S^{n-1}$, where
$K^* = F^{\rm st}(K,S^0)$ is the $S$-dual of $K$.
For spectra $A$ and $B$,  let $F(A,B)$ denote the
associated function spectrum. Then $\pi_0(F(A,B)) = [A,B]$.

The first step is to rewrite 
$$
F(E,D_2(E)) \simeq  F(E,E\smsh E)_{h\Bbb Z_2} 
$$
(the 
$\Bbb Z_2$-action on $F(E,E\smsh E)$ is induced by permutation
action on the smash product $E \smsh E$.) To see this,
note there is a natural map from right to left.
That this map is a  weak equivalence can established
by induction on a cell structure for $E$, recalling that
$E$ is homotopy finite.

Substituting in the value of $E$ into the above,
we get
$$
F(E,D_2(E)) \simeq 
F(K^*\smsh S^{n-1},(K^*\smsh S^{n-1})^{\smsh 2})_{h\Bbb Z_2} \, .
$$
Now, using the fact that $S^{n-1} \smsh S^{n-1}$ with permutation action
is homeomorphic to $S^{(n{-}1)\alpha}\smsh S^{n{-}1}$ with diagonal
action, the right side of the last display can be rewritten as 
$$
 F(K^*, S^{(n{-}1)\alpha} \smsh K^* \smsh K^*)_{h\Bbb Z_2} \, .
$$
For homotopy finite spectra $A$ and $B$, it is well known that
the transpose map $F(A,B) \to F(B^*,A^*)$ is a weak
equivalence. Consequently, there
is a $\Bbb Z_2$-equivariant weak equivalence of spectra
$$
F^{\rm st}(K\smsh K,S^{(n{-}1)\alpha} \smsh K) \,\, \simeq \,\, 
  F(K^*, S^{(n{-}1)\alpha} \smsh K^* \smsh K^*) \, .
$$
given by the transpose map.

Taking homotopy orbits of this last equivalence,
and assembling the prior information,
we conclude that there is a weak equivalence
of spectra
$$
F(E,D_2(E)) \simeq F^{\rm st}(K\smsh K,S^{(n{-}1)\alpha} \smsh K)_{h\Bbb Z_2} \, .
$$
Applying $\pi_0$ to this last equivalence completes the proof.
\end{proof}

To complete the proof
of Theorem \ref{enumerate}, one just needs to apply
Corollary \ref{reformulation}, Lemma \ref{first}, Lemma \ref{second}
and Theorem \ref{desuspend} in the stated order
(to apply \ref{desuspend}, use
the fact that $E = F^{\rm st}(K,S^{n{-}1})$ is $(n{-}k{-}2)$-connected
and $\dim E \le n - r - 2$). We leave it to the reader to check
that the inequalities listed in the statement of Theorem \ref{enumerate}
suffice to apply these results.

\section{Proof of Theorem \ref{exist}}

The proof of Theorem \ref{exist} is almost identical
to the proof of  Theorem \ref{enumerate}. There are two essential
differences: the first is that 
instead of using Theorem \ref{desuspend}, we need to use
the following existence result for realizing a spectrum as a
suspension spectrum in the metastable range:

\begin{thm}[Klein \cite{Klein_susp}] 
There is 
an obstruction
$$
\delta_E \in [E,\Sigma D_2(E)]\, ,
$$
(depending only on the homotopy type of $E$) which 
is trivial whenever $E$ has the homotopy type of a suspension spectrum.

Conversely, if $E$ is $r$-connected, $r \ge 1$ and $\dim E \le 3r{+}2$, then
$E$ has the homotopy type of a suspension spectrum if $\delta_E = 0$.
\end{thm}

The second essential difference is that when 
$E = F^{\rm st}(K,S^{n{-}1})$,
we have an isomorphism of abelian groups
$$
[E,\Sigma D_2(E)] \cong 
F^{\rm st}(K\smsh K,S^{(n{-}1)\alpha + 1}\smsh K)_{h\Bbb Z_2}\, .
$$
The obstruction $\theta_K$ is defined so as to correspond to the obstruction
$\delta_E$
with respect to this isomorphism of abelian groups.
We omit the details.

\section{Applications to two cell complexes}

\subsection*{Existence}
It seems that case of embedding complexes with two cells was first considered
by Cooke \cite{Cooke_embed1} (see also \cite{Cooke_embed2}) and later 
by Connolly and Williams \cite[\S5]{Connolly-Williams}.

Let $K = S^p \cup_f e^{q+1}$ be a two cell complex, where
$f\: S^q \to S^p$ is some map. Let $E := F^{\rm st}(K,S^{n{-}1})
$ denote the stable Spanier-Whitehead 
$(n-1)$-dual of $K$. Set $p' = n {-}p{-}2$ and $q' = n{-}q{-}2$.

Then $E$ is the homotopy cofiber of a stable 
umkehr map
$$
f^*\:S^{p'} \to S^{q'}\, .
$$
As stable classes in $\pi_{q{-}p}^{\rm st}(S^0)$, we have
$$
[f^*] = [\pm f]\, .
$$ 
Tracing through the definition of the umkehr map, with
slightly extra care, 
the sign can be determined as $(-1)^{qp'}$.

In any case,
$E$ has the homotopy type of a suspension spectrum if
and only if $f^*$ is represented by an
unstable map. In our range, this is equivalent to
demanding that the  James-Hopf invariant 
$$
H_2(f^*) =     \pi_{p'}^{\rm st}(D_2(S^{q'}))
$$
is trivial.

\subsection*{Enumeration}
Suppose $K = S^p \cup_f e^{q{+}1}$ admits an embedding
in $S^n$. An analysis similar to the previous case shows that there is
an isomorphism of based sets
$$
E(K,S^n) \cong \pi_{p'+1}^{\rm st}(D_2(S^{q'}))
$$
At the prime 2, the stable homotopy groups appearing
on the right have been calculated
by Mahowald in degrees $p' \le \min(3q'{-}3,2q' {+}29)$
(see Mahowald \cite[table 4.1]{Mahowald_memoir}).

For example, suppose that $q' \equiv 1 \text{ mod } 16$.
Then the  first few groups are
\medskip
$$
\vbox{
\offinterlineskip \tabskip=2pt
\halign{
\strut # &
\vrule # &
\hfil # \hfil &
\vrule # &
\hfil # \hfil &
\vrule # &
\hfil # \hfil &
\vrule # &
\hfil # \hfil &
\vrule # &
\hfil # \hfil &
\vrule # &
\hfil # \hfil &
\vrule # &
\hfil # \hfil &
\vrule # &
\hfil # \hfil &
\vrule # &
\hfil # \hfil &
\vrule #\cr
\omit& \multispan{17}{\hrulefill}\cr
& & $j$ & & $0$
& & 1& &  2 & & 3 && 4 && 5 && 6&\cr
\omit& \multispan{17}{\hrulefill}\cr
& & $\pi_{2q'{+} j}(D_2(S^{q'}))$ &&  $\Bbb Z_2$ & & $\Bbb Z_2$
& & $\Bbb Z_8$ & &  $\Bbb Z_2$& & $0$ && $\Bbb Z_2$
&& $\Bbb Z_{16} \oplus \Bbb Z_2 $
&\cr
\omit& \multispan{17}{\hrulefill}\cr
}}\, .
$$ 
\medskip

\section{Embeddings of Poincar\'e spaces}

In this section we assume that
$K$ is a $r$-connected Poincar\'e duality space
of formal dimension $k$.
\medskip

\begin{rems_no}  The
Browder-Casson-Sullivan-Wall theorem (\cite[Th.\ 12.1]{Wall_2nd})
says that concordance classes of Poincar\'e embeddings of $K$ in $S^n$ 
are in one to one correspondence with embeddings up to homotopy of
$K$ in $S^n$.

If $K$ is a closed PL manifold, then
\cite[Th.\ 11.3.1]{Wall_2nd} implies that
$E(K,S^n)$ is in bijection with the isotopy classes of 
locally flat PL embeddings of $K$ in $S^n$.
\end{rems_no}
\medskip

By \cite[Lem.\ 2.8]{Wall_2nd}, we can find a homotopy
equivalence $K \simeq L \cup e^k$,
where $L$ is a finite complex and $\dim L \le k {-} r {-} 1$.
In particular, we have a cofibration sequence of $\Bbb Z_2$-spaces
$$
L \smsh K \cup_{L \smsh L} K \smsh L \to K \smsh K \to S^k \smsh S^k \, .
$$
The first term of this sequence has dimension $\le 2k {-}r {-}1$,
so we may infer that the evident map 
$$
F^{\rm st}(S^k \smsh S^k,S^{(n{-}1)\alpha {+} 1}\smsh K)_{h\Bbb Z_2}
\to 
F^{\rm st}(K \smsh K,S^{(n{-}1)\alpha {+} 1}\smsh K)_{h\Bbb Z_2}
$$
is $(n{-}2(k{-}r){+}1)$-connected. In particular, if $n \ge 2(k{-}r)$,
we see that this map induces an isomorphism on path components.

By elementary manipulations which we omit, 
there is an identification 
$$
F^{\rm st}(S^k \smsh S^k,S^{(n{-}1)\alpha {+} 1}\smsh K)_{h\Bbb Z_2}
\,\, \simeq \,\,
F^{\rm st}(S^{n-2}, K \smsh D_2(S^{n{-}k{-}1})) \, .
$$
We conclude:

\begin{thm} \label{pd-exist}
Assume in addition $n \ge 2(k{-}r)$.
Then the obstruction
$\theta_K$  is detected in the abelian group
$$
\pi_{n-2}^{\rm st}(K \smsh D_2(S^{n{-}k{-}1}))  \, .
$$
\end{thm}
\medskip

\begin{rem_no} Let $\nu$ be the Spivak
normal fibration of $K$; we consider $\nu$ has having fiber
a {\it stable} $(-k)$-sphere. 
Let $K^\nu$ denote the Thom spectrum of $\nu$. 
When $K$ embeds in $S^n$, the fibration $\nu$
compresses down to an
{\it unstable} $(n{-}k{-}1)$-spherical fibration. Conversely, when
$\nu$ compresses, a construction due to Browder gives
an  embedding of $K$ in $S^{n{+}1}$ (see \cite{Browder_embed}).

It is therefore tempting  to try and relate $\theta_K$ 
to the obstruction theoretic problem of finding a
compression of $\nu$. We do not as yet have a solution to this.
\end{rem_no} 
\bigskip

By essentially the same argument that
proves \ref{pd-exist}, we have
 
\begin{thm}  \label{pd-enumerate} Assume $n > 2(k{-}r)$.
Then
the function $\phi_K$ can be rewritten as 
$$
\phi_K\:E(K,S^n) \to 
\pi_{n{-}1}^{\rm st}(K \smsh D_2(S^{n{-}k{-}1})) \, .
$$
\end{thm}
\medskip

The remainder of this section is devoted to obtaining
corollaries of \ref{pd-exist} and \ref{pd-enumerate}.
Our first result shows that $\phi_K$
is homological in the fringe dimension beyond the
stable range.

\begin{cor}[Compare {\cite[Th.\ 2.3]{Haef-Hirsch}}, \cite{Habegger_embed}]
\label{Haef-Hirsch} The obstruction $\phi_K$ to embedding
$K$ in $S^{2k{-}r{-}1}$ lives in the abelian group
$$
H_{r{+}1}(K;\Bbb Z_s)\, ,
$$
where $s = 1 + (-1)^{k{-}r{+}1}$. 
\end{cor}

\begin{proof}The Hurewicz map 
\begin{align*}
\pi_{2k{-}r{-}3}^{\rm st}(K \smsh D_2(S^{k{-}r{-}2}))
& \to \,\,
H_{2k{-}r{-}3}(K \smsh D_2(S^{k{-}r{-}2})) \\
& \cong \,\, 
H_{r{+}1}(K) \otimes H_{2(k{-}r{-}2)}(D_2(S^{k{-}r{-}2}))  \\
& \cong \,\, H_{r{+}1}(K;\Bbb Z_s)
\end{align*}
is an isomorphism in this degree. Now apply \ref{pd-exist}.
\end{proof}
\medskip

By a similar argument 
which we omit (use \ref{pd-enumerate}), we obtain

\begin{cor} [Compare {\cite[Th.\ 2.4]{Haef-Hirsch}}, \cite{Habegger_embed}]
The set of concordance classes of embeddings of $K$
in $S^{2k{-}r{+}2}$ is isomorphic to
$$
H_{r{+}1}(K;\Bbb Z_s) \,  ,
$$
where $s = 1 + (-1)^{k{-}r}$.
\end{cor}
\medskip

Our next pair of corollaries concern the outcome 
of tensoring with the rationals.

\begin{cor} \label{rational-existence} 
If $n \equiv k \text{ \rm mod 2}$, then 
$\theta_K\otimes \Bbb Q$ is trivial.
Otherwise, 
$\theta_K \otimes \Bbb Q$ is 
detected in the vector space $H_{2k{-}n}(K;\Bbb Q)$.
\end{cor}
\medskip

\begin{proof} 
If $n \equiv k \text{ \rm mod } 2$, then 
$\pi_*(D_2(S^{n-k-1}))\otimes \Bbb Q$ 
is trivial.  We infer that $\pi_*(K\smsh D_2(S^{n-k-1}))\otimes \Bbb Q$
is also trivial.
The first part now follows using \ref{pd-exist}.

For the second part, note that the transfer
$$D_2(S^{n{-}k{-}1}) \to (S^{n{-}k{-}1})^{\smsh 2}$$
is, rationally, the inclusion of a wedge summand. Smashing with $K$ 
and applying rational homotopy, we
infer that 
$\pi_{n{-}2}^{\rm st}(K \smsh D_2(S^{n{-}k{-}1}))\otimes \Bbb Q$ is
a summand of 
$\pi_{n{-}2}^{\rm st}(K \smsh (S^{n{-}k{-}1})^{\smsh 2})\otimes \Bbb Q$. 
Over the rationals, stable homotopy coincides with homology.
It follows that $\theta_K\otimes \Bbb Q$ 
is detected in $H_{2k{-}n}(K;\Bbb Q)$.
\end{proof}

\begin{cor} \label{rational-enumeration} Assume $K$ embeds in $S^n$.
Assume inequality \eqref{number2} holds. Then
$E(K,S^n)$ is finitely generated.

If $n \equiv k \text{ \rm mod } 2$, then 
 $E(K,S^n)$ is finite. Otherwise,
$E(K,S^n)\otimes \Bbb Q$  is a direct summand of 
$H_{2k{-}n{+}1}(K;\Bbb Q)$.
\end{cor}

\begin{proof}[Proof of \ref{rational-enumeration}] 
The first part follows from \ref{pd-enumerate} because 
$\pi_{n{-}1}^{\rm st}(K \smsh D_2(S^{n{-}k{-}1}))$
is finitely generated. 
The second part is proved in a manner similar to
\ref{rational-existence}. We omit the details.
\end{proof}

A direct consequence of \ref{rational-enumeration}
is:

\begin{cor} \label{pd_embed_finiteness} Assume the inequality \eqref{number2} holds. If
the betti number $b_{2k{-}n{+1}}(K)$ is trivial, then
there are finitely many concordance classes of embeddings
of $K$ in $S^n$.
\end{cor}

This last result gives Corollary \ref{isotopy} of the introduction
using the remarks about manifolds given at the beginning of this section.

\section{Localization at 2}
Let $K$ and $K'$ be $r$-connected finite complexes
with  $\dim K,\dim K' \le k$.

\begin{thm} \label{2-local}
Suppose that $f\: K \to K'$ is a 2-local homotopy
equivalence. Assume that inequality \eqref{number1} holds.
Then $K$ embeds in $S^n$ if and only if $K'$ does.
\end{thm}
\medskip

\begin{rem_no} Rigdon \cite{Rigdon_2local} and
Williams \cite{Williams_2local} prove 
a similar result for manifolds in the metastable range $n \ge 3/2(k{+}1)$.
The main difference between their result and ours is 
that ours holds outside of the metastable range
at the expense of an additional connectivity hypothesis.
\end{rem_no} 
\medskip

\begin{proof}[Proof of \ref{2-local}] The induced map of stable $(n{-}1)$-duals
$$
E':= F^{\rm st}(K',S^{n{-}1})\overset{f^*}\to F^{\rm st}(K,S^{n{-}1}) =: E
$$
is clearly a 2-local equivalence.
By \cite[Th.\ D]{Klein_susp}, $E'$ is a suspension spectrum if and only if
$E$ is. The result now 
follows by applying lemmas \ref{first}, \ref{sw=d} and 
Theorem \ref{frank-bruce}.
\end{proof}
\bigskip

%mention manifolds

%rationalize or invert two 
%?? K --> K' 2 local equiv then K embeds <=> K' does 
% can we simplify the connectivity hypotheses?

%\begin{thebibliography}
%\bibliographystyle{invent}
%\bibliography{john}
%\end{thebibliography}
\end{document}